\documentclass[11pt]{article}

\usepackage{url}
\usepackage{fullpage}
\usepackage[english]{babel}
\usepackage{amsmath}
\usepackage{amsfonts}
\usepackage{amsthm}
\usepackage{times}

\newcommand{\Z}{\mathbb Z}
\newcommand{\Q}{\mathbb Q}
\newcommand{\R}{\mathbb R}
\newcommand{\GL}{\mathsf{GL}}
\newcommand{\SC}{{\mathbf\Delta}}
\newcommand{\MT}{{\mathcal T}}

\newcommand{\gldz}{\GL_d(\Z)}
\newcommand{\gldr}{\GL_d(\R)}

\newcommand{\seight}{{\mathcal S}^8}
\newcommand{\sd}{{\mathcal S}^d}
\newcommand{\sdo}{{\mathcal S}^d_{>0}}
\newcommand{\sdgeo}{{\mathcal S}^d_{\geq 0}}

\renewcommand{\vec}[1]{\boldsymbol{#1}}

\DeclareMathOperator{\trace}{trace}
\DeclareMathOperator{\vol}{vol}
\DeclareMathOperator{\conv}{conv}
\DeclareMathOperator{\Co}{Co}
\DeclareMathOperator{\Del}{Del}
\DeclareMathOperator{\diag}{diag}

\theoremstyle{definition}
\newtheorem{definition}{Definition}[section]
\newtheorem{proposition}[definition]{Proposition}
\newtheorem{lemma}[definition]{Lemma}
\newtheorem{theorem}[definition]{Theorem}
\newtheorem{corollary}[definition]{Corollary}

\title{Local Covering Optimality of Lattices:\\
Leech Lattice versus Root Lattice $\mathsf{E}_8$}

\author{Achill Sch\"urmann, Frank Vallentin
\footnote{Partially supported by the Edmund Landau Center
for Research in Mathematical Analysis and Related Areas, sponsored by
the Minerva Foundation (Germany).}}
        
\date{10th November 2004}

\begin{document}

\maketitle

\begin{abstract}
We show that the Leech lattice gives a sphere covering which is
locally least dense among lattice coverings.  We show that a similar
result is false for the root lattice $\mathsf{E}_8$. For this we
construct a less dense covering lattice whose Delone subdivision has a
common refinement with the Delone subdivision of $\mathsf{E}_8$. The
new lattice yields a sphere covering which is more than $12\%$ less
dense than the formerly best known given by the lattice
$\mathsf{A}_8^*$.  Currently, the Leech lattice is the first and only
known example of a locally optimal lattice covering having a
non-simplicial Delone subdivision.  We hereby in particular answer a
question of Dickson posed in 1968.  By showing that the Leech lattice
is rigid our answer is even strongest possible in a sense.
\end{abstract}

\section{Introduction}
\label{sec:intro}

The Leech lattice is the exceptional lattice in dimension $24$. Soon
after its discovery by Leech~\cite{leech-1967} it was conjectured that
it is extremal for several geometric problems in $\R^{24}$: the
kissing number problem, the sphere packing problem and the sphere
covering problem.

In 1979, Odlyzko and Sloane and independently Levenshtein solved the
kissing number problem in dimension $24$ by showing that the Leech
lattice gives an optimal solution. Two years later, Bannai and Sloane
showed that it gives the unique solution up to isometries (see
\cite{cs-1988}, Ch.\ 13, 14). Unlike the kissing number problem, the
other two problems are still open.

Recently, Cohn and Kumar~\cite{ck-2004} showed that the Leech lattice
gives the unique densest \textit{lattice} sphere packing in $\R^{24}$.
Furthermore they showed, that the density of any sphere packing
(without restriction to lattices) in $\R^{24}$ cannot exceed the one
given by the Leech lattice by a factor of more than $1+1.65 \cdot
10^{-30}$.

At the moment it is not clear how one can prove a corresponding result
for the sphere covering problem. In this paper we take a first step
into this direction by showing

\begin{theorem}
\label{th:leech}
The Leech lattice gives a sphere covering which is locally optimal
among lattices.
\end{theorem}

In Section~\ref{sec:basics} we give precise definitions of
``locally optimality'' and of all other terms needed.

Surprisingly, a result similar to Theorem~\ref{th:leech} does not hold
for the root lattice $\mathsf{E}_8$, which is the exceptional lattice
in dimension $8$. As the Leech lattice, $\mathsf{E}_8$ gives the
unique solution to the kissing number problem (see~\cite{cs-1988},
Ch.\ 13, 14) and is conjectured to be extremal for the sphere packing
problem.  Blichfeldt~\cite{blichfeldt-1934} showed that it gives an
optimal lattice sphere packing. Later, Vetchinkin~\cite{vetchinkin-1980}
showed that sphere packings of all other $8$-dimensional lattices are
less dense.  Besides giving another proof of this, Cohn and Kumar
\cite{ck-2004} demonstrated that the density of a sphere packing in
$\R^8$ cannot exceed the one of $\mathsf{E}_8$ by a factor of more
than $1+10^{-14}$.  In contrast to the Leech lattice, $\mathsf{E}_8$
cannot be an extremal sphere covering in its dimension though. Conway
and Sloane note in \cite{cs-1988}, Ch.\ 2: "It is surprising that
$\mathsf{A}^*_8$, with [covering density] $\Theta = 3.6658\ldots$, is
better than $\mathsf{E}_8$, which has $\Theta = 4.0587\ldots$".  Here,
we show that $\mathsf{E}_8$ does not even give a locally optimal
lattice sphere covering, by constructing a new $8$-dimensional
lattice, which yields a less dense sphere covering and whose Delone
subdivision has a common refinement with the one of $\mathsf{E}_8$.

\begin{theorem}
\label{th:newcovering}
There exists a lattice with covering density $\Theta < 3.2013$ whose
Delone subdivision has a common refinement with the one of
$\mathsf{E}_8$.
\end{theorem}

Note that this new sphere covering beats the old record holder
$\mathsf{A}^*_8$ in dimension $8$.  The proof of
Theorem~\ref{th:newcovering} relies on computational methods we
developed in \cite{sv-2004} and up to now we can only give an
approximation of the best ``known'' covering lattice.  In any case, by
Proposition~\ref{prop:local-min} due to Barnes and Dickson
Theorem~\ref{th:newcovering} yields:

\begin{corollary}
The root lattice $\mathsf{E}_8$ does not give a locally optimal
lattice sphere covering.
\end{corollary}

By Theorem~\ref{th:leech}, we give a first example of a locally optimal
covering lattice whose Delone subdivision is not simplicial. By this
we give an affirmative answer to a question of Dickson posed in
\cite{dickson-1968}. The Leech lattice gives even a strongest possible
example, in the sense that it is rigid (see Section~\ref{sec:basics}
for details).  Our proof in Section~\ref{sec:rigid} immediately
applies to $\mathsf{E}_8$, giving a new proof of $\mathsf{E}_8$'s
rigidity, first observed by Baranovskii and Grishukhin~\cite{bg-2002}.

\begin{theorem}
\label{th:rigid}
The Leech lattice and $\mathsf{E}_8$ are rigid. 
\end{theorem}

\section{Lattices, Positive Quadratic Forms and Delone Subdivisions}
\label{sec:basics}

In this section we briefly review some concepts and results about
lattices and their relation to positive definite quadratic forms (PQFs
from now on). For further reading we refer to \cite{cs-1988},
Ch.~2~\S2.2 and \cite{sv-2004}.

Let $\R^d$ denote a $d$-dimensional Euclidean space with unit ball
$B^d$. Its volume is $\kappa_d={\pi^{d/2}}/{\Gamma(d/2+1)}$. A (full
rank) \textit{lattice} $L$ is a discrete subgroup in $\R^d$, that is,
there exists a regular matrix $A \in \gldr$ with $L = A\Z^d$.  The
\textit{determinant} $\det(L) = |\det(A)|$ of $L$ is independent of
the chosen basis~$A$.  The Minkowski sum $L + \alpha B^d = \{\vec{v} +
\alpha\vec{x} : \vec{v} \in L, \vec{x} \in B^d\}$, $\alpha \in
\R_{>0}$, is called a \textit{lattice packing} if the translates of
$\alpha B^d$ have mutually disjoint interiors and a \textit{lattice
covering} if $\R^d = L + \alpha B^d$.  The \textit{packing radius}
\[
\lambda(L)
=
\max\{ \lambda : \text{$L + \lambda B^d$ is a lattice packing} \},
\]
and the \textit{covering radius}
\[
\mu(L)
=
\min \{ \mu : \text{$L + \mu B^d$ is a lattice covering} \},
\]
are both attained. They are homogeneous and therefore, the \textit{covering
density} $\Theta(L)=\frac{\mu(L)^d\kappa_d}{\det(L)}$ and the \textit{packing
density} $\delta(L)=\frac{\lambda(L)^d\kappa_d}{\det(L)}$ are invariant
with respect to scaling of~$L$.

Given a $d$-dimensional lattice $L=A\Z^d$ with basis $A$ we associate
a $d$-dimensional PQF $ Q[\vec{x}] =\vec{x}^t A^t A \vec{x}
=\vec{x}^t G \vec{x} $, where the \textit{Gram matrix} $G = A^t A$ is
symmetric and positive definite.  We will carelessly identify
quadratic forms with symmetric matrices by saying $Q = G$ and
$Q[\vec{x}] = \vec{x}^t Q \vec{x}$.  The set of quadratic forms is a
${d+1 \choose 2}$-dimensional real vector space $\sd$, in which the
set of PQFs forms an open, convex cone~$\sdo$.  The PQF~$Q$ depends on
the chosen basis $A$ of $L$. For two arbitrary bases $A$ and $B$ of
$L$ there exists a $U\in\gldz$ with $A = BU$.  Thus, $\gldz$ acts on
$\sdo$ by $ Q \mapsto U^t Q U$. A PQF Q can be associated to different
lattices $L = A\Z^d$ and $L' = A'\Z^d$. In this case there exists an
orthogonal transformation $O$ with $A = OA'$. Note that the packing
and covering density are invariant with respect to orthogonal
transformations.

The \textit{determinant} (or discriminant) of a PQF $Q$ is defined by
$\det(Q)$.  The \textit{homogeneous minimum} $\lambda(Q)$ and the
\textit{inhomogeneous minimum} $\mu(Q)$ are given by
\[
\lambda(Q) = \min_{\vec{v} \in \Z^d \setminus \{\vec{0}\}} Q[\vec{v}],
\qquad
\mu(Q) = \max_{\vec{x} \in \R^d} \min_{\vec{v} \in \Z^d} Q[\vec{x} - \vec{v}].
\]
If $Q$ is associated to $L$, then $\det(L) = \sqrt{\det(Q)}$,
$\mu(L) = \sqrt{\mu(Q)}$, $\lambda(L) = \sqrt{\lambda(Q)}/2$.

We say that a lattice $L$ with associated PQF $Q$ gives a
\textit{locally optimal lattice covering} or a \textit{locally optimal
lattice packing}, if there is a neighborhood of $Q$ in $\sdo$, so
that for all $Q'$ in the neighborhood we have $\Theta(Q) \leq
\Theta(Q')$, respectively $\delta(Q) \geq \delta(Q')$.

A polytope $P = \conv\{\vec{v}_1, \ldots, \vec{v}_n\}$, with
$\vec{v}_1, \ldots, \vec{v}_n \in \Z^d$, is called a \textit{Delone
polytope} of~$Q$ if there exists a $\vec{c} \in \R^d$ and a real
number $r \in \R$ with $Q[\vec{v}_i - \vec{c}] = r^2$ for all $i = 1,
\ldots, n$, and for all other lattice points $\vec{v} \in \Z^d
\setminus \{\vec{v}_1, \ldots, \vec{v}_n\}$ we have strict inequality
$Q[\vec{v} - \vec{c}] > r^2$. The set of all Delone polytopes is
called the \textit{Delone subdivision} of~$Q$. Note that the
inhomogeneous minimum of $Q$ is at the same time the maximum squared
circumradius of its Delone polytopes. We say that the Delone
subdivision of a PQF $Q'$ is a \textit{refinement} of the Delone
subdivision of $Q$, if every Delone polytope of $Q'$ is contained in a
Delone polytope of $Q$.

By a theory of Voronoi~\cite{voronoi-1908} (see also \cite{sv-2004}),
the set of PQFs with a fixed Delone subdivision is an open polyhedral
cone in $\sdo$ --- the \textit{secondary cone} of the subdivision.
In the literature the secondary cone is sometimes called
$L$-type domain of the subdivision.  The topological closures of
these secondary cones give a face-to-face tessellation of $\sdgeo$,
the set of all positive semi-definite quadratic forms.  The relative
interior of a face in this tessellation contains PQFs that have the
same Delone subdivision. If a face is contained in the boundary of a
second face, then the corresponding Delone subdivision of the first is
a true refinement of the second one. The relative interior of faces of
minimal dimension $1$ contain \textit{rigid} PQFs. They have the
special property that every PQF $Q'$ in a sufficiently small
neighborhood and not being a multiple of $Q$, has a Delone subdivision
which is a true refinement of $Q$'s subdivision.  The relative
interior of faces of maximal dimension ${d+1 \choose 2}$ contain PQFs
whose Delone subdivision is a triangulation, that is, it consists of
simplices only. We refer to such a subdivision as a \textit{simplicial} 
Delone subdivision or \textit{Delone triangulation}.

We transfer the terminology of Delone subdivisions from PQFs to
lattices by saying that the Delone subdivision of the lattice $L'$ is
a \textit{refinement} of the Delone subdivision of the lattice $L$, if
there are associated PQFs $Q'$ and $Q$ so that the Delone subdivision
of $Q'$ is a refinement of the Delone subdivision of $Q$.  A lattice
is called \textit{rigid} if an associated PQF is rigid.

In \cite{dickson-1968}, Dickson states: ``Whether it is possible for
$f_0$ [a PQF giving a locally optimal lattice covering] to occur on
the boundary of a cone [...] is still a matter of conjecture.''  Hence
to answer his question affirmatively, one has to find a PQF giving a
locally optimal lattice covering with non-simplicial Delone
subdivision.  The following proposition by Barnes and Dickson (see
\cite{bd-1967}, \cite{dickson-1968} \S5) shows that there is
essentially at most one such lattice covering for each Delone
subdivision:

\begin{proposition} 
\label{prop:local-min}
A lattice $L$ (or an associated PQF) gives a locally optimal lattice
sphere covering iff it minimizes the covering density among all
lattices whose Delone subdivisions have a common refinement with $L$'s
Delone subdivision.  Moreover, such a lattice $L$ is determined
uniquely up to dilations and orthogonal transformations.
\end{proposition}

\section{The Leech Lattice, the Root Lattice $\mathsf{E}_8$, and their Rigidity}
\label{sec:rigid}

Let us introduce the two lattices. We gathered the information mostly
from \cite{cs-1988}, Ch.~4~\S8,\S11.

The Leech lattice $\Lambda$ with associated PQF $Q_{\Lambda}$
satisfies $\det(Q_{\Lambda})=1$, $\lambda(Q_{\Lambda})=4$ and
$\mu(Q_{\Lambda})=2$.  Thus, its packing density, already given by
Leech \cite{leech-1967}, is $\delta(Q_{\Lambda}) = \kappa_{24}$ and
best possible among all lattices in $\R^{24}$ (\cite{ck-2004},
Th.~9.3).  Its covering density is $\Theta(Q_{\Lambda}) = 4096 \cdot
\kappa_{24}$. The first proof of this fact is due to Conway, Parker
and Sloane \cite{cs-1988}, Ch.~23. There they also classified the $23$
different (up to congruences) Delone polytopes of $Q_{\Lambda}$
attaining the maximum circumradius $\sqrt{2}$.

From their list, we will consider Delone polytopes of type
$\mathsf{A}^{24}_1$ to prove the rigidity in this section and those of
type $\mathsf{A}_{24}$ to prove the local optimality in
Section~\ref{sec:leech}. To describe them we define $Q_{\Lambda}(n) =
\{\vec{v} \in \Z^{24} : Q_{\Lambda}[\vec{v}] = 2n\}$. The Delone
polytopes of type $\mathsf{A}^{24}_1$ are $24$-dimensional regular
cross polytopes with respect to the metric induced by $Q_{\Lambda}$.
They are of the form $\vec{v} + \conv\{\vec{v}_0, \ldots,
\vec{v}_{47}\}$, $\vec{v} \in \Z^{24}$, where $\vec{v}_0 = 0$,
$\vec{v}_{24} \in Q_{\Lambda}(4)$ and all the other $\vec{v}_i \in
Q_{\Lambda}(2)$ satisfy $\vec{v}_j + \vec{v}_{j+24} = \vec{v}_{24}$,
$j = 0,\ldots, 47$ (indices computed modulo $48$). The Delone
polytopes of type $\mathsf{A}_{24}$ are $24$-dimensional simplices
having $275$ edge vectors in $Q_{\Lambda}(2)$ and $25$ in
$Q_{\Lambda}(3)$. This is the only information about $\mathsf{A}_{24}$
we will need in Section~\ref{sec:leech}.

For the proof of rigidity it is convenient to work with the following
coordinates with respect to the standard basis of $\R^{24}$:
The vectors of squared length $4$ of $\Lambda$ are
of shape $\frac{1}{\sqrt{8}}((\pm 4)^20^{22})$,
$\frac{1}{\sqrt{8}}((\pm 2)^80^{16})$ and
$\frac{1}{\sqrt{8}}(\mp3(\pm 1)^{23})$, where permitted permutations of
coordinates and permitted positions of minus signs are explained in
\cite{cs-1988}, Ch.~10 ~\S3.2. Here, we only need the $2^2\binom{24}{2}$
vectors of the first type, where all permutations of coordinates
and all positions of minus signs are allowed.

The root lattice $\mathsf{E}_8$ with associated PQF $Q_{\mathsf{E}_8}$
satisfies $\det(Q_{\mathsf{E}_8}) = 1$, $\lambda(Q_{\mathsf{E}_8}) =
2$ and $\mu(Q_{\mathsf{E}_8}) = 1$. Thus, its packing density is
$\delta(Q_{\mathsf{E}_8}) = \frac{1}{16} \cdot \kappa_{8}$ and best
possible among all lattices in $\R^{8}$ (\cite{blichfeldt-1934}).  Its
covering density is $\Theta(Q_{\mathsf{E}_8}) = \kappa_{8}$. The first
proof of this fact is due to Coxeter~\cite{coxeter-1946}, who gave a
complete description of the Delone subdivision of
$Q_{\mathsf{E}_8}$. It is a tiling of $\R^8$ into regular simplices
and regular cross polytopes with respect to the metric induced by
$Q_{\mathsf{E}_8}$.  Define $Q_{\mathsf{E}_8}(n) = \{\vec{v} \in
\Z^{8} : Q_{\mathsf{E}_8}[\vec{v}] = n\}$. Note that this slightly
differs from the definition of $Q_{\Lambda}(n)$ above.  Then the cross
polytopes are of the form $\vec{v} + \conv\{\vec{v}_0, \ldots,
\vec{v}_{15}\}$, $\vec{v} \in \Z^{8}$, where $\vec{v}_0 = 0$,
$\vec{v}_{8} \in Q_{\mathsf{E}_8}(4)$ and all the other $\vec{v}_i \in
Q_{\mathsf{E}_8}(2)$ satisfy $\vec{v}_j + \vec{v}_{j+8} =
\vec{v}_{8}$, $j = 0,\ldots, 15$ (indices computed modulo $16$). The
regular simplices have $36$ edge vectors in $Q_{\mathsf{E}_8}(2)$. We
shall give more information about this Delone subdivision in
Section~\ref{sec:newcovering}.

For the proof of $\mathsf{E}_8$'s rigidity and for the construction of
a new sphere covering in Section~\ref{sec:newcovering} it is again
convenient to work with explicit coordinates with respect to the
standard basis of $\R^8$. Set
\begin{equation}
\label{eq:e8}
\mathsf{E}_8 = \{\vec{x} \in \R^8 : 
\text{$\vec{x} \in \Z^8 \cup (\frac{1}{2} + \Z)^8$ and  
$\sum_{i=1}^8 x_i \in 2\Z$}\}.
\end{equation}
The automorphism group of $\mathsf{E}_8$ is generated by all
permutations of the $8$ coordinates, by all even sign changes and by
the matrix $H = \diag(H_4,H_4)$ where
\[
H_4 = \frac{1}{2}
\left(
\begin{array}{rrrr}
1 & 1 & 1 & 1 \\
1 &-1 & 1 &-1 \\
1 & 1 &-1 &-1 \\
1 &-1 &-1 & 1 \\
\end{array}
\right).
\]
There are $240$ vectors of squared length $2$ in $\mathsf{E}_8$:
$2^2\binom{8}{2}$ of shape $((\pm 1)^20^6)$ and $2^7$ of shape $((\pm
\frac{1}{2})^8)$ where the number of minus signs is even. The $2160$
vectors of squared length $4$ are: $2 \cdot 8$ of shape $((\pm 2)0^7)$,
$2^4\binom{8}{4}$ of shape $((\pm 1)^40^4)$ and $2^7 \cdot 8$ of shape
$(\pm \frac{3}{2} (\pm \frac{1}{2})^7)$ where the number of minus
signs is odd.

Now we proceed to the proof of Theorem~\ref{th:rigid}.  We will handle
both cases simultaneously. For this, denote by $L$ the Leech lattice
$\Lambda$ or the root lattice $\mathsf{E}_8$ and write $d$ for the
rank of $L$. We shall show that every PQF $Q$, whose Delone
subdivision contains the above mentioned cross polytopes, is a
multiple of $Q_{L}$. By $\langle \cdot, \cdot\rangle$ we denote the
inner product given by $Q_{L}$, i.e.\ $\langle \vec{x}, \vec{y}
\rangle = \vec{x}^t Q_{L} \vec{y}$, and by $(\cdot, \cdot)$ we denote
the inner product given by $Q$.

Let $\vec{v}, \vec{w} \in Q_{L}(2)$ with $\langle \vec{v}, \vec{w}
\rangle = 0$.  So, $\vec{v} + \vec{w} \in Q_{L}(4)$. Therefore, 
$\vec{0}, \vec{v} + \vec{w}, \vec{v}, \vec{w}$ are vertices of
a Delone cross polytope, as considered above.  Let $\vec{c}$ be the
center of its circumsphere, hence $Q[\vec{0} - \vec{c}] = Q[\vec{v} +
\vec{w} - \vec{c}] = Q[\vec{v} - \vec{c}] = Q[\vec{w} - \vec{c}]$. 
Then a straightforward calculation reveals $(\vec{v}, \vec{w}) = 0$.

Now we switch to coordinates with respect to the standard basis.  
Choosing a basis $A$ of $L$ gives the associated PQF $Q_L = A^t A$.
Obviously, $Q = A^t (A^t)^{-1} Q A^{-1} A$. We denote the entries 
of $C = (A^t)^{-1} Q A^{-1}$ by $(c_{ij})$ and by $\vec{e}_i$ we denote 
the $i$-th canonical basis vector of $\R^d$. Let $\vec{v}_i \in \R^d$ 
be the coordinate vector of $\vec{e}_i$ with respect to the basis $A$,
that is $A\vec{v}_i = \vec{e}_i$. Then by our choice of coordinates 
and by the argument above we have for $i \neq j$ the orthogonality
\[
0 = 
(\vec{e}_i + \vec{e}_j)^t (\vec{e}_i - \vec{e}_j) =
(\vec{v}_i + \vec{v}_j)^t Q_L (\vec{v}_i - \vec{v}_j) =
(\vec{v}_i + \vec{v}_j)^t Q (\vec{v}_i - \vec{v}_j) =
(\vec{e}_i + \vec{e}_j)^t C (\vec{e}_i - \vec{e}_j) =
c_{ii} - c_{jj}.
\]
Moreover, for pairwise different indices
$i,j,k,l$, the orthogonality yields
\[
0 = (\pm \vec{e}_i \pm \vec{e}_j)^t (\pm \vec{e}_k \pm \vec{e}_l) =
\pm c_{ik} \pm c_{il} \pm c_{jk} \pm c_{jl}.\] 
Hence, the matrix $C$
is a multiple of the identity matrix and so $Q$ is a multiple of
$Q_L$, which proves Theorem~\ref{th:rigid}.

\section{Local Lower Bounds for the Covering Density}
\label{sec:locallowerbound}

In this section we briefly describe a variant of a method due to
Ryshkov and Delone which enables us to compute local lower bounds for
the covering density. This is a slight variation of the method
described in \cite{sv-2004} and \cite{vallentin-2003}.

Let $L = \conv\{\vec{v}_1, \ldots, \vec{v}_{d+1}\} \subseteq \R^d$ be
a simplex in the $d$-dimensional Euclidean space with inner product
given by the PQF $Q$. Its \textit{centroid} is $\vec{m} =
\frac{1}{d+1}\sum_i \vec{v}_i$. Let $\vec{c}$ be the center of its
circumsphere and let $r$ be its circumradius. Using Apollonius'
formula (see \cite{berger-1987} \S 9.7.6) we get
\[
r^2 = Q[\vec{c} - \vec{m}] + \frac{1}{(d+1)^2} \sum_{k \neq l}
Q[\vec{v}_k - \vec{v}_l].
\] 

\begin{proposition}
\label{prop:mubound}
Let $L_1 = \conv\{\vec{v}_{1,1}, \ldots, \vec{v}_{1,d+1}\}, \ldots,
L_n = \conv\{\vec{v}_{n,1}, \ldots, \vec{v}_{n,d+1}\}$ be a collection
of Delone simplices of $Q$ with radii $r_1, \ldots, r_n$. Then, 
the inhomogeneous minimum is bounded by
\begin{equation}
\label{eq:mu}
\mu(Q) \geq \max\limits_{i} r_i^2 \geq
\frac{1}{n(d+1)^2}\sum_i\sum_{k \neq l} Q[\vec{v}_{i,k} -
\vec{v}_{i,l}].
\end{equation}
\end{proposition}

The proof is straightforward.  We can use the foregoing proposition to
get local lower bounds for the covering density of PQFs having $L_1,
\ldots, L_n$ as Delone simplices. We fix the determinant of $Q$ and
minimize the right hand side of (\ref{eq:mu}), which is a linear
function:

\begin{proposition}
\label{prop:optimize_equidiscriminant}
Let $D>0$. A linear function $f(Q') = \trace(FQ')$ with a PQF $F$ has
a unique minimum on the determinant $D$ surface $\{Q \in \sdo:
\text{$\det(Q) = D$}\}$. Its value is $d \sqrt[d]{D \det F}$ and the
minimum is attained by the PQF $\sqrt[d]{D \det F}F^{-1}$.
\end{proposition}

\noindent For a proof of
Proposition~\ref{prop:optimize_equidiscriminant} we refer to
\cite{vallentin-2003}, Proposition~8.2.2. Together,
Proposition~\ref{prop:mubound} and
Proposition~\ref{prop:optimize_equidiscriminant} yield

\begin{corollary}
\label{cor:moments}
As in Proposition~\ref{prop:mubound}, let $L_1, \ldots, L_n$ be a
collection of Delone simplices of a PQF $Q$. Then 
\[
\Theta(Q) \geq \sqrt{\left(\frac{d}{d+1}\right)^d \det F} \cdot
\kappa_d,
\] 
with $F=\frac{1}{n(d+1)} \sum_i\sum_{k \neq l} (\vec{v}_{i,k} -
\vec{v}_{i,l})(\vec{v}_{i,k} - \vec{v}_{i,l})^t$ which is a PQF.
\end{corollary}

\section{Local Optimality of the Leech Lattice}
\label{sec:leech}

For the proof of Theorem~\ref{th:leech} we use the fact that any
non-empty set $\frac{1}{\sqrt{2n}}Q_{\Lambda}(n)$, $n > 0$, forms a
spherical $11$-design (\cite{cs-1988}, Ch.~7, Th.~23) in the
Euclidean space with inner product $\langle \cdot, \cdot \rangle$.
Generally, a \textit{spherical $t$-design} $X$ is a non-empty finite
subset of the unit sphere $S^{d-1}=\{\vec{x}\in\R^d :
\langle\vec{x},\vec{x}\rangle=1\}$ satisfying $\frac{1}{\vol
S^{d-1}}\int_{S^{d-1}} f(\vec{x}) d\vec{x} =
\frac{1}{|X|}\sum_{\vec{x} \in X} f(\vec{x})$ for every polynomial $f
: \R^d \to \R$ of degree at most $t$. Here, $\vol S^{d-1}$ denotes the
surface volume of $S^{d-1}$, not the volume of the enclosed
ball. Equivalently, $X$ is a spherical $t$-design iff it satisfies
the equalities (see \cite{venkov-2001}, Th.~3.2):
\[
\begin{array}{ll}
\sum_{\vec{x} \in X} \langle \vec{x},\vec{y} \rangle^k = 0, &
\text{for all $\vec{y} \in \R^d$ and all odd $k \leq t$,}\\
\sum_{\vec{x} \in X} \langle \vec{x},\vec{y} \rangle^k = 
\frac{1\cdot3\cdots(k-1)}{d(d+2)\cdots(d+k-2)} |X| \langle \vec{y},\vec{y}\rangle^{k/2}, &
\text{for all $\vec{y} \in \R^d$ and all even $k \leq t$.}
\end{array}
\]
For the proof of Theorem~\ref{th:leech} the following spherical
$2$-design property is even sufficient:

\begin{lemma}
\label{lemma:venkov}
Let $Q\in\sdo$ and let $X\subset\R^d$ denote a spherical $2$-design
with respect to the inner product given by $Q$. Then
\[
\sum_{\vec{x} \in X} \vec{x}\vec{x}^t
= \frac{|X|}{d} Q^{-1}
.
\]
\end{lemma}

\begin{proof}
Since $X$ forms a spherical $2$-design, we have
$\sum_{\vec{x} \in X} \left(\vec{x}^t Q \vec{y}\right)^2=
\frac{|X|}{d} (\vec{y}^t Q \vec{y})$.
On the other hand 
\[
\sum_{\vec{x} \in X} \left(\vec{x}^t Q \vec{y}\right)^2 =
\sum_{\vec{x} \in X} \vec{y}^t Q (\vec{x}\vec{x}^t) Q \vec{y} =
\vec{y}^t Q \left( \sum_{\vec{x} \in X} \vec{x}\vec{x}^t \right) Q
\vec{y} .
\]
Thus because both identities are valid for all $\vec{y}\in\R^d$ we
derive the equality stated in the lemma.
\end{proof}

Now we finish the proof of Theorem~\ref{th:leech}. Let $L$ be a Delone
simplex of $Q_{\Lambda}$ of type $\mathsf{A}_{24}$. We apply
Corollary~\ref{cor:moments} to the orbit of $L$ under the automorphism
group $\Co_0 = \{T \in \GL_{24}(\Z) : T^t Q_{\Lambda} T =
Q_{\Lambda}\}$ of $Q_{\Lambda}$.

For every PQF $Q$ for which the simplices $TL$, $T \in \Co_0$, are
Delone simplices, we have $\Theta(Q) \geq \sqrt{(\frac{24}{25})^{24}
\det F} \cdot \kappa_{24}$ with $F = \frac{1}{25|\Co_0|} \sum_{T \in
\Co_0} \sum_{\vec{e}} \vec{e}\vec{e}^t$, where $\vec{e}$ runs through
all the edge vectors of $TL$.  Since $L$ has $275$ edges in
$Q_{\Lambda}(2)$ and $25$ edges in $Q_{\Lambda}(3)$ and because of the
transitivity of $\Co_0$ on $Q_{\Lambda}(2)$ and $Q_{\Lambda}(3)$
(\cite{cs-1988}, Ch.~10, Th.~27) we get
$$ F = \frac{1}{25|\Co_0|} \left( \frac{275 |\Co_0|}{|Q_{\Lambda}(2)|}
\sum_{\vec{e} \in Q_{\Lambda}(2)} \vec{e}\vec{e}^t + \frac{25
|\Co_0|}{|Q_{\Lambda}(3)|} \sum_{\vec{e} \in Q_{\Lambda}(3)}
\vec{e}\vec{e}^t \right).
$$ By Lemma~\ref{lemma:venkov} (applied to $Q_{\Lambda}/4$ and
$Q_{\Lambda}/6$) this yields
$$ F = \frac{1}{25} \left(\frac{275}{|Q_{\Lambda}(2)|} \cdot
\frac{|Q_{\Lambda}(2)|}{6} Q_{\Lambda}^{-1} + \frac{25}{|Q_{\Lambda}(3)|} \cdot
\frac{|Q_{\Lambda}(3)|}{4} Q_{\Lambda}^{-1}\right) = \frac{5^2}{2^2\cdot3}
Q_{\Lambda}^{-1}
$$ Since $\det Q_{\Lambda}^{-1} = 1$, it follows $\det F =
\frac{5^{48}}{2^{48}3^{24}}$ and finally, by
Corollary~\ref{cor:moments}, we derive $\Theta(Q) \geq 4096 \cdot
\kappa_{24} = \Theta(Q_{\Lambda})$.

\section{A New Sphere Covering in Dimension 8}
\label{sec:newcovering}

If we apply the method of Theorem~\ref{th:leech} to $\mathsf{E}_8$ we
get a local lower bound of $\sqrt{(8/9)^8} \cdot \kappa_8 \approx
0.6243 \cdot \kappa_8$. But $\Theta(\mathsf{E}_8) = \kappa_8$, since
the circumradius of the regular cross polytopes in $\mathsf{E}_8$'s
Delone subdivision is~$1$. Despite this gap, $\mathsf{E}_8$ could be a
locally optimal covering lattice. The following proof of
Theorem~\ref{th:newcovering} shows that this is not the case.  By
Proposition~\ref{prop:local-min} we have to find a PQF $Q$ with
$\Theta(Q) < \Theta(Q_{\mathsf{E}_8})$ so that $\Del(Q)$ and
$\Del(Q_{{\mathsf E}_8})$ have a common refinement.

Below, we describe a systematic way to attain a refining Delone
triangulation of $\Del(Q_{{\mathsf E}_8})$. Given such a triangulation
we can find an approximation of the unique PQF minimizing $\Theta$
among all PQFs in the closure of its secondary cone by solving a
convex programming problem on a computer. We give such an
approximation in Appendix~\ref{sec:verification} and verify its
properties with a simple computer program. Since we carry out the
verification using exact arithmetic only, the proof of
Theorem~\ref{th:newcovering} is rigorous.

First, we describe all $\Z^8$-periodic triangulations refining the
Delone subdivision of $Q_{{\mathsf E}_8}$, that is, all sets
$\mathcal{P}$ of simplices satisfying the following conditions:
\begin{description}
\item[no additional vertices:] all vertices of simplices $L\in\mathcal{P}$ lie in $\Z^8$. 
\item[periodicity:] $\forall L \in \mathcal{P}, \vec{v} \in \Z^8: \vec{v} + L \in \mathcal{P}$.
\item[face-to-face tiling:] $\forall L, L' \in \mathcal{P}: L \cap L' \in \mathcal{P}$.
\item[refinement:] $\forall L \in \mathcal{P} \; \exists L' \in \Del(Q_{\mathsf{E}_8}): L \subseteq L'$.
\item[covering:] $\forall \vec{x} \in \R^8 \; \exists L \in \mathcal{P}: \vec{x} \in L$.
\end{description}
Recall from Section~\ref{sec:rigid} that $\Del(Q_{\mathsf{E}_8})$
consists of simplices and cross polytopes only. Thus for a
$\Z^8$-periodic triangulation refining $\Del(Q_{{\mathsf E}_8})$ we
have to specify how to split the $8$-dimensional cross polytopes into
simplices.

We say that two polytopes $P$ and $P'$ are \textit{$\Z^8$-equivalent}
if $P' = \vec{v} + P$ for some $\vec{v} \in \Z^8$.  Every $\vec{w} \in
Q_{{\mathsf E}_8}(4)$ defines a Delone cross polytope
$P_{\vec{w}}=\conv\{\vec{v}_0,\ldots,\vec{v}_{15}\}$ with
$\vec{v}_0=\vec{0}$, $\vec{v}_8=\vec{w}$ and all other $\vec{v}_j\in
Q_{{\mathsf E}_8}(2)$ with $\vec{v}_j+\vec{v}_{j+8}=\vec{v}_8$
(indices computed modulo $16$).

Two $\vec{w}, \vec{w}' \in Q_{{\mathsf E}_8}(4)$ define
$\Z^8$-equivalent cross polytopes iff $\vec{w}' \in \vec{w} + 2\Z^8$,
because then the difference $\frac{1}{2}\vec{w}'-\frac{1}{2}\vec{w}$
of their centers is in $\Z^8$.  Under this equivalence relation the
set $Q_{{\mathsf E}_8}(4)$ splits into $135$ classes, containing $8$
pairs of mutually orthogonal vectors $\pm \vec{w}_1,\dots,\pm
\vec{w}_8$. Each of the $\vec{w}_i$ equivalent to $\vec{w}$ is a
diagonal of $P_{\vec{w}}$, e.g.~$\vec{w}_i=\vec{v}_i-\vec{v}_{i+8}$,
$i=0,\dots,7$. In the coordinate system introduced in~(\ref{eq:e8}) the $135$
classes are (see \cite{cs-1988}, Ch.~6, \S3):
\begin{description}
\item[] $1$ class: $\pm 2\vec{e}_1, \ldots, \pm 2\vec{e}_8$.
\item[] $70$ classes: $8$ elements $\pm\vec{e}_a \pm \vec{e}_b \pm
\vec{e}_c \pm \vec{e}_d$ with an even number of minus signs and $8$
elements $\pm \vec{e}_e \pm \vec{e}_f \pm \vec{e}_g \pm \vec{e}_h$
with an even number of minus signs and with $\{a,\ldots, h\} =
\{1,\ldots, 8\} $; or the same with an odd number of minus signs.
\item[] $64$ classes: $8$ pairs of vectors of shape
$(\pm \frac{3}{2} (\pm \frac{1}{2})^7)$ with odd number of minus signs,
where the position of $\pm\frac{3}{2}$ is permuted to all $8$ 
coordinates.
\end{description}       
We can split each cross polytopes $P_{\vec{w}}$ into simplices in
eight different ways by adding a diagonal.  Without loss of generality
we add the diagonal $\conv\{\vec{v}_0, \vec{v}_{8}\}$ and split the
cross polytope $P_{\vec{w}}$ into the $128$ simplices $\conv\{\vec{0},
\vec{v}_{8}, \vec{v}_{j_1}, \ldots, \vec{v}_{j_7}\}$, where ${j_k}\in
\{k,k+8\}$.  Thus altogether we get $8^{135}$ different
$\Z^8$-periodic triangulations refining $\Del(Q_{{\mathsf E}_8})$.

Now, which of these periodic triangulations are Delone triangulation
for some PQF?  To decide this, we take a closer look at the tiling
$\Del(Q_{\mathsf{E}_8})$ and at secondary cones $\SC(\MT)$ of Delone
triangulations $\MT$ refining $\Del(Q_{\mathsf{E}_8})$.

We already described the cross polytopes of
$\Del(Q_{\mathsf{E}_8})$. Centers of simplices of
$\Del(Q_{\mathsf{E}_8})$ containing the origin are the $17280$ vectors
$\frac{1}{3}\vec{v}$, where $\vec{v}$ is a vector of
$Q_{\mathsf{E}_8}(8)$ not in $2Q_{\mathsf{E}_8}(2)$.  We say two
polytopes are adjacent in the tiling, if they share a facet. Each
simplex is adjacent to $9$ cross polytopes and each cross polytope is
adjacent to $128$ simplices and $128$ cross polytopes. A simplex and a
cross polytope both containing the origin are adjacent iff the inner
product, with respect to $Q_{\mathsf{E}_8}$, of their centers equals
$\frac{5}{6}$.  Two cross polytopes both containing the origin are
adjacent iff the inner product of their centers equals $\frac{3}{4}$.

For some computations it is useful to have coordinates of vertices,
with respect to the coordinate system (\ref{eq:e8}): The vertices of the
cross polytope $P$ defined by the center $\vec{e}_1$ are $\vec{0},
2\vec{e}_1, \vec{e}_1 \pm \vec{e}_i$, $i = 2, \ldots, 8$. An adjacent
Delone cross polytope is defined by the center $\vec{c} =
(\frac{3}{4}, -\frac{1}{4}, \frac{1}{4}, \frac{1}{4}, \frac{1}{4},
\frac{1}{4}, \frac{1}{4}, \frac{1}{4})$. Its vertices are $\vec{0},
2\vec{c}, \vec{e}_1 - \vec{e}_2, 2\vec{c} - (\vec{e}_1 - \vec{e}_2),
\vec{e}_1 + \vec{e}_i, 2\vec{c} - (\vec{e}_1 + \vec{e}_i)$, $i =
3,\ldots,8$. A Delone simplex adjacent to $P$ is defined by the center
$\vec{c}' =
(\frac{5}{6},\frac{1}{6},\frac{1}{6},\frac{1}{6},\frac{1}{6},\frac{1}{6},\frac{1}{6},\frac{1}{6})$. Its
vertices are $\vec{0}, (\frac{1}{2},\ldots,\frac{1}{2}), \vec{e}_1 +
\vec{e}_i$, $i = 2,\ldots,8$. Since the automorphism group of
$\mathsf{E}_8$ acts transitively on vectors of squared length $4$ and
since the stabilizer of $\pm2\vec{e}_1$ in $\mathsf{E}_8$'s
automorphism group is the group generated by even sign changes and by
permutations of the last $7$ coordinates, the knowledge of the
coordinates given is enough to describe the whole Delone subdivision.

The secondary cones $\SC(\MT)$ are open polyhedral cones and by the
theory of Voronoi, they are given by linear forms on $\seight$ called
\textit{regulators}.  Each pair of adjacent simplices
$L=\conv\{\vec{v}_0,\ldots,\vec{v}_8\}$,
$L'=\conv\{\vec{v}_1,\ldots,\vec{v}_{9}\}$ gives a regulator
$\varrho_{(L,L')}$. If $\alpha_0,\ldots,\alpha_{9} \in \Q$ are the
uniquely determined numbers with $\alpha_0=1$, $\sum_{i=0}^{9}\alpha_i
= 0$ and $\sum_{i=0}^{9}\alpha_i \vec{v}_i = \vec{0}$, then $
\varrho_{(L,L')}(Q) = \sum_{i=0}^{9}\alpha_i Q[\vec{v}_i] $ for
$Q\in\seight$ and
\[
\SC(\MT)
=\left\{
Q\in\seight :
\text{$\varrho_{(L,L')}(Q) > 0$, $(L,L')$ pair of adjacent simplices of $\MT$}
\right\}.
\]
Note that $\varrho_{(L+v,L'+v)}=\varrho_{(L,L')}$ for all $\vec{v} \in
\Z^8$.

For a triangulation $\MT$ which is a refinement of
$\Del(Q_{\mathsf{E}_8})$ we distinguish between three types of pairs
of adjacent simplices $(L,L')$. In the first case one of the simplices
is a simplex of $\Del(Q_{\mathsf{E}_8})$ and the other one is not. In
the two other cases both simplices are not simplices of
$\Del(Q_{\mathsf{E}_8})$. In the second case they refine adjacent
cross polytopes, in the third case they refine the same cross
polytope. In the first two cases we have
$\varrho_{(L,L')}(Q_{\mathsf{E}_8})>0$ and in the last case
$\varrho_{(L,L')}(Q_{\mathsf{E}_8})=0$. Since $\mathsf{E}_8$ is rigid,
the closures of secondary cones of Delone triangulations refining
$\Del(Q_{\mathsf{E}_8})$ cover a sufficient small neighborhood of the
ray containing multiples of $Q_{\mathsf{E}_8}$. Thus, $\MT$ is a
Delone triangulation for some PQF refining $\Del(Q_{\mathsf{E}_8})$
iff
\[
\left\{ Q\in\seight \; :\; \varrho_{(L,L')}(Q)>0, \; (L,L')
\;\text{refining the same cross polytope of }\Del(Q_{\mathsf{E}_8})
\right\}
\]
is not empty.

The three types of regulators are easily computed, e.g.\ with help of
the coordinates given above. In the first case, let $P =
\conv\{\vec{v}_0,\ldots,\vec{v}_{15}\}$ be a cross polytope of
$\Del(Q_{{\mathsf E}_8})$ with the notational convention: $\vec{v}_8
\in Q_{\mathsf{E}_8}(4)$, $\vec{v}_i + \vec{v}_{i+8} = \vec{v}_8$. Let
$L' = \conv\{\vec{v}'_0, \ldots, \vec{v}'_8\}$ be a simplex of
$\Del(Q_{{\mathsf E}_8})$ with $\vec{v}'_i = \vec{v}_i$, $i = 0,
\ldots, 7$. Let $\vec{c}$ be the centroid of $P$ and $\vec{c}'$ be the
centroid of $L'$. Then $\vec{c}'=\frac{1}{9}\left( \vec{v}'_0+\ldots + 
\vec{v}'_8 \right)$ and $ \frac{1}{4} \vec{c} +
\frac{3}{4} \vec{c}' = \frac{1}{8}(\vec{v}_0 + \cdots +
\vec{v}_7)$. Suppose the edge $\conv\{\vec{v}_k, \vec{v}_{k+8}\}$, $k
\in \{0,\ldots,7\}$, belongs to $\MT$. We have $\vec{c} =
\frac{1}{2}(\vec{v}_k + \vec{v}_{k+8})$ and we derive $
\frac{1}{8}(\vec{v}_k + \vec{v}_{k+8}) + \frac{1}{12}(\vec{v}'_8 +
\vec{v}_0 + \cdots + \vec{v}_7) = \frac{1}{8}(\vec{v}_0 + \cdots +
\vec{v}_7)$. Therefore we get the regulator
\begin{equation}
\label{eq:regulator1}
\varrho_{(L,L')}(Q) = Q[\vec{v}_k] + Q[\vec{v}_{k+8}] + \frac{2}{3} Q[\vec{v}'_8] -
\frac{1}{3}Q[\vec{v}_0] - \cdots - \frac{1}{3}Q[\vec{v}_7].
\end{equation}

In the second case, let $P_1=\conv\{\vec{v}_0,\ldots,\vec{v}_{15}\}$
and $P_2=\conv\{\vec{v}'_0,\ldots,\vec{v}'_{15}\}$ be two adjacent
cross polytopes with the usual notational convention and with
$\vec{v}_i=\vec{v}'_i$ for $i=0,\ldots,7$. Then the centers
$c=\frac{1}{2}\vec{v}_8$ and $c'=\frac{1}{2}\vec{v}'_8$ of the cross
polytopes satisfy the relation $\frac{1}{2}\left(c+c'\right) =
\frac{1}{8}\left(\vec{v}_0+\cdots+\vec{v}_7\right)$. Let us assume
that the diagonals $\conv\{\vec{v}_{k},\vec{v}_{k+8}\}$ and
$\conv\{\vec{v}'_{k'},\vec{v}'_{k'+8}\}$ with $k,k'\in\{0,\dots,7\}$
belong to $\MT$. Then, since $c$, $c'$ are the centers of these
diagonals, we derive $\vec{v}_k + \vec{v}_{k+8} + \vec{v}'_k +
\vec{v}'_{k'+8} = \frac{1}{2}\vec{v}_0+\cdots+\frac{1}{2}\vec{v}_7 $.
Therefore we get the regulator
\begin{equation}
\label{eq:regulator2}
\varrho_{(L,L')}(Q) = Q[\vec{v}_{k}] + Q[\vec{v}_{k+8}] +  Q[\vec{v}'_{k'}] + Q[\vec{v}'_{k'+8}]
-\frac{1}{2}Q[\vec{v}_{0}] - \cdots - \frac{1}{2}Q[\vec{v}_{7}]                      
.
\end{equation}

In the third case, let $P = \conv\{\vec{v}_0,\ldots,\vec{v}_{15}\}$ be a
cross polytope with the usual notational convention.  Then adjacent
simplices are of the form $L=\conv\{\vec{0}, \vec{v}_{8},
\vec{v}_{j_1}, \ldots, \vec{v}_{j_7}\}$, $L'=\conv\{\vec{0},
\vec{v}_{8}, \vec{v}_{j'_1}, \ldots, \vec{v}_{j'_7}\}$, where
$j_k,j'_k\in \{k,k+8\}$ and $j'_k=j_k+8$ only for one
$k\in\{1,\ldots,7\}$.  Because of
$\vec{v}_{j_k}+\vec{v}_{j_k+8}=\vec{v}_0 + \vec{v}_8$ we get seven
regulators
\begin{equation}
\label{eq:regulator3}
\varrho_{(L,L')}(Q) = Q[\vec{v}_{j_k}] + Q[\vec{v}_{j_k+8}] -Q[\vec{v}_0] - Q[\vec{v}_8]
,\quad k=1,\dots,7.
\end{equation}
Note that these conditions are equivalent to $Q[\vec{v}_8]<
Q[\vec{v}_{j_k}-\vec{v}_{j_{k+8}}]$, $k=1,\dots,7$, which means that
the chosen diagonal $\vec{v}_8$ is shorter than the other seven with
respect to the metric induced by $Q$.

We tried to generate all Delone triangulations refining
$\Del(Q_{\mathsf{E}_8})$ by an exhaustive computer search. But this
seems to be hopeless since they are far to many. So we decided to
generate a Delone triangulation which has a fairly large symmetry
group. For this we choose a subgroup $G$ of $Q_{\mathsf{E}_8}$'s
automorphism group which, in the coordinate system~(\ref{eq:e8}), is
generated by permutations of the last $7$ coordinates and by the
involution $\vec{x} \mapsto -\vec{x}$.

\begin{proposition}
\label{prop:periodic-triang}
There are exactly four $\Z^8$-periodic triangulations refining
$\Del(Q_{\mathsf{E}_8})$ invariant under the group $G$. Exactly
two of them are Delone triangulations and both are equivalent under
the action of $Q_{\mathsf{E}_8}$'s full automorphism group.
\end{proposition}

\begin{proof}
To show that there is essentially one Delone triangulation with the
prescribed symmetries, we will again work with the
coordinate system~(\ref{eq:e8}). In Table~1 we list the orbits of squared
length $4$ vectors under the action of $G$.
\begin{center}
\begin{tabular}[t]{|c|c|c|}
\hline
\textbf{\#} & \textbf{representative} & \textbf{orbit size}\\
\hline
1. & $2000000$ & $1$ \\
2. & $0200000$ & $7$ \\
3. & $11110000$ & $\binom{7}{3}$ \\
4. & $1\overline{1}110000$ & $7\binom{6}{2}$ \\
5. & $1\overline{11}10000$ & $7\binom{6}{2}$ \\
6. & $1\overline{111}0000$ & $\binom{7}{3}$\\
\hline
\end{tabular}
\hspace{2ex}
\begin{tabular}[t]{|c|c|c|}
\hline
7. & $01111000$ & $\binom{7}{4}$ \\
8. & $0\overline{1}111000$ & $7 \binom{6}{3}$ \\
9. & $0\overline{11}11000$ & $\frac{1}{2}\binom{7}{2} \binom{5}{2}$ \\
10. & $\frac{1}{2}(3\overline{1}111111)$ & $\binom{7}{1}$ \\
11. & $\frac{1}{2}(3\overline{111}1111)$ & $\binom{7}{3}$ \\
12. & $\frac{1}{2}(3\overline{11111}11)$ & $\binom{7}{5}$ \\
13. & $\frac{1}{2}(3\overline{1111111})$ & $\binom{7}{7}$\\
\hline
\end{tabular}
\hspace{2ex}
\begin{tabular}[t]{|c|c|c|}
\hline
14. & $\frac{1}{2}(13\overline{1}11111)$ & $7\binom{6}{1}$\\
15. & $\frac{1}{2}(13\overline{111}111)$ & $7\binom{6}{3}$ \\
16. & $\frac{1}{2}(13\overline{11111}1)$ & $7\binom{6}{5}$ \\
17. & $\frac{1}{2}(1\overline{3}111111)$ & $7\binom{6}{0}$ \\
18. & $\frac{1}{2}(1\overline{311}1111)$ & $7\binom{6}{2}$ \\
19. & $\frac{1}{2}(1\overline{31111}11)$ & $7\binom{6}{4}$ \\
20. & $\frac{1}{2}(1\overline{3111111})$ & $7\binom{6}{6}$ \\
\hline
\end{tabular}\\
\vspace*{0.2cm}
\textsf{\textbf{Table 1.}} 
Orbits of squared length $4$ vectors in $\mathsf{E}_8$.
Minus signs are given by bars.
\end{center}

\noindent
To define a $\Z^8$-periodic triangulation refining
$\Del(Q_{\mathsf{E}_8})$ we have to choose a collection of orbits
$(O_i)_{i \in I}$, $I \subseteq \{1,\ldots, 20\}$, so that for every
of the $135$ classes of possible diagonals $C$ we have
$\left|\bigcup_{i \in I} O_i \cap C \right| = 2$. This restriction
immediately gives $\{2, 4, 5, 8, 9, 14, 15, 16, 17, 18, 19\} \cap I =
\emptyset$. For example, the four vectors $\pm 02000000$, $\pm
00200000$ are in $O_2 \cap C$. On the other hand we have to have $\{1,
6, 11, 12, 13\} \subseteq I$. Now there are two binary choices left:
either we have $3 \in I$ or $7 \in I$ and either we have $10 \in I$ or
$20 \in I$.  From these four triangulations only those given by $I_1 =
\{1, 3, 6, 10, 11, 12, 13\}$ and $I_2 = \{1, 6, 7, 11, 12, 13, 20\}$
are Delone triangulations: Under the prescribed symmetry we can assume
that there are numbers $\alpha$, $\beta$, $\gamma$, $\delta$ with
\[
\alpha = (\vec{e}_1, \vec{e}_1), \;\;
\beta = (\vec{e}_1, \vec{e}_i),\;\;
\gamma = (\vec{e}_i, \vec{e}_i),\;\;
\delta = (\vec{e}_i, \vec{e}_j),\;\; i = 2,\ldots, 8, j = i + 1, \ldots, 8.
\]
Suppose we choose orbit $3$. By (\ref{eq:regulator3}) this implies the
inequality
\[
(\frac{1}{2}(11111111),\frac{1}{2}(1111\overline{1111})) = \frac{1}{4}(\alpha + 6\beta - \gamma - 6\delta) < 0,
\]
then choosing orbit $20$ implies
\[
\begin{array}{l}
(\frac{1}{2}(1\overline{1}000000), \frac{1}{2}(\overline{11111111})) = \frac{1}{4}(-\alpha - 6\beta + \gamma + 6\delta) < 0,
\end{array}
\]
yielding a contradiction. Hence we have to choose orbit $10$
instead. A similar calculation shows that if we choose orbit $7$, then
we have to choose orbit $20$. To see that these triangulations are
Delone triangulations we still have to give a PQF satisfying all
regulators in (\ref{eq:regulator3}). We postpone this to 
Appendix~\ref{sec:verification}.

By applying the transformation $HAH$, where $A$ exchanges the first
and fifth coordinate and their signs, and $H$ is the transformation $H
= \diag(H_4, H_4)$ we see that both Delone triangulations are
equivalent. The transformation $HAH$ is an element of $\mathsf{E}_8$'s
automorphism group exchanging the relevant orbits of vectors of size
$2$ by $O_1 \leftrightarrow O_{13}$, $O_3 \leftrightarrow O_7$, $O_6
\leftrightarrow O_{11}$, $O_{10} \leftrightarrow O_{20}$, $O_{12}
\leftrightarrow O_{12}$.
\end{proof}

Given the Delone triangulation $\MT$ refining
$\Del(Q_{\mathsf{E}_8})$, attained in this way or another, we can
compute an approximation of the unique PQF minimizing the covering
density $\Theta$ among all PQFs in the closure of its secondary
cone. For details we refer to \cite{sv-2004} and give only a brief
sketch. We have to solve the optimization problem: maximize $\det Q$
where $Q$ lies in the closure of the secondary cone $\SC(\MT)$ and the
circumradius of every simplex in $\MT$ with respect to $Q$ is bounded
by $1$. This is a convex optimization problem and we can approximate
the solution using a computer. With help of the software
\texttt{MAXDET} written by Wu, Vandenberghe and Boyd (see
\cite{vbw-1998}) we found a PQF $\tilde{Q}$ with covering density
$\Theta \approx 3.2012$. But \texttt{MAXDET} uses floating point
arithmetic. So we have to verify that the Delone subdivision of the
new PQF has a common refinement with $\Del(Q_{\mathsf{E}_8})$ and we
have to verify $\tilde{Q}$'s covering density. We did this by writing
a program which uses only rational arithmetic. We give more details
and the PQF in Appendix~\ref{sec:verification}. The successful
verification proves in particular Theorem~\ref{th:newcovering}.

Finally, we report on some numerical evidences. The PQF $\tilde{Q}$
lies on the boundary of the secondary cone $\SC(\MT)$. The closure of
$\SC(\MT)$ has $428$ facets and only one of these facets does not
contain $Q_{\mathsf{E}_8}$. We applied the optimization to the Delone
triangulation which belongs to the secondary cone adjacent to this
facet. There we found a PQF with covering density $\Theta \approx
3.1423$, which is the best known covering density in dimension $8$ so
far.  It seems that the PQF we approximated by $\tilde{Q}$ is not
locally optimal. This would follow by Proposition
\ref{prop:local-min}, if we knew that the by $\tilde{Q}$ approximated
PQF lies on the boundary of $\SC(\MT)$. In any case we are left with
the open problem to find a globally best covering lattice in dimension
$8$.  Currently, we do not know where to search for such a lattice.

\section{Remarks on the Packing-Covering Problem}
\label{sec:remarks}

Together with the local optimality of the Leech lattice with respect
to the packing problem, we immediately derive the local optimality of
the Leech lattice with respect to another problem (see
\cite{sv-2004}): The lattice packing-covering problem asks to
minimize the packing-covering constant
$\gamma(L)=\mu(L)/\lambda(L)=(\Theta(L)/\delta(L))^{1/d}$ among all
$d$-dimensional lattices $L$. For the Leech lattice we derive
$\gamma(\Lambda)=\sqrt{2}$ which again is a strict local minimum. We
can offer two different proofs for this: either by using the global
optimality of the Leech lattice with respect to the packing problem or
by deriving a local lower bound along the lines of the proof of
Theorem~\ref{th:leech}.  Here, analogous tools to Proposition
\ref{prop:optimize_equidiscriminant} and Corollary~\ref{cor:moments}
(see \cite{sv-2004}, Sec.~10) are needed.

For $\mathsf{E}_8$ the situation seems to differ from the covering
case: Given a fixed Delone triangulation, the problem of finding the
minimum $\gamma$ among all lattices with the same Delone triangulation
can also be formulated as a convex optimization problem (see
\cite{sv-2004}).  In contrast to the covering case, the application of
\texttt{MAXDET} to the triangulation constructed in
Section~\ref{sec:newcovering} indicates that
$\gamma(\mathsf{E}_8)=\sqrt{2}$ may in fact be a local optimum, as
conjectured by Zong \cite{zong-2002}. He conjectured that
$\mathsf{E}_8$ gives the global optimum.

Note the remarkable fact that $d=1,2$ are the only known cases where
the minimum covering density and the maximum packing density --- and
therefore the minimum packing-covering constant
$\gamma_d=\min_L\gamma(L)$ --- are known to be attained by the same
lattice. Maybe yet another exceptional property of the beautiful Leech
lattice\dots

\section*{Acknowledgments}

We like to thank Mathieu Dutour and the anonymous referee for valuable
comments on a previous version.

\begin{appendix}

\section{Verification of Numerical Results}
\label{sec:verification}

The PQF $\tilde{Q}$ of Section~\ref{sec:newcovering} is $\tilde{Q}=34229189769 Q_1 - 17121746137 Q_2$
with
\[
\begin{array}{rcl}
Q_1 & = &
\begin{pmatrix}
    1 &    0 &    0 &    0 &    0 &    0 &    0 &    2\\
    0 &  4/7 & -2/3 &    0 &    0 &    0 &    0 &    0\\
    0 & -2/3 &  4/3 & -2/3 &    0 &    0 &    0 &    0\\
    0 &    0 & -2/3 &  4/3 & -2/3 &    0 &    0 &    0\\
    0 &    0 &    0 & -2/3 &  4/3 & -2/3 &    0 &    0\\
    0 &    0 &    0 &    0 & -2/3 &  4/3 & -2/3 &    0\\
    0 &    0 &    0 &    0 &    0 & -2/3 &  4/3 &    0\\
    2 &    0 &    0 &    0 &    0 &    0 &    0 &    4\\
\end{pmatrix}
\quad\mbox{and}\\
Q_2 & = &
\begin{pmatrix}
    0 &    1 &    0 &    0 &    0 &    0 &    0 &  7/2\\
    1 &    0 & -2/3 &    0 &    0 &    0 &    0 &    0\\
    0 & -2/3 &  4/3 & -2/3 &    0 &    0 &    0 &    0\\
    0 &    0 & -2/3 &  4/3 & -2/3 &    0 &    0 &    0\\
    0 &    0 &    0 & -2/3 &  4/3 & -2/3 &    0 &    0\\
    0 &    0 &    0 &    0 & -2/3 &  4/3 & -2/3 &    0\\
    0 &    0 &    0 &    0 &    0 & -2/3 &  4/3 &    0\\
  7/2 &    0 &    0 &    0 &    0 &    0 &    0 &    7\\
\end{pmatrix}.
\end{array}
\]
Here, $Q_1$ and $Q_2$ form a basis of the subspace of $\seight$
invariant under the group~$G$ intersected with the subspace given by
regulators $\varrho_{(L,L')}$ of (\ref{eq:regulator1}),
(\ref{eq:regulator2}), (\ref{eq:regulator3}) with
$\varrho_{(L,L')}(\tilde{Q}) = 0$.  To give a rigorous proof of
Theorem~\ref{th:newcovering} we have to verify that one of the two triangulations
we constructed in Section~\ref{sec:newcovering} is a refining Delone
triangulation of $\tilde{Q}$'s Delone subdivision and that its
covering density is at most $3.2013$.  To do this we supply the
\texttt{MAGMA} program \texttt{newcovering8.m} available from the
\texttt{arXiv.org} e-print archive. To access it, download the source
files for the paper \texttt{math.MG/0405441}. There we also included
the Gram matrix of the ``best known'' covering lattice in the
additional file \texttt{currentbest8.txt}.

The program \texttt{newcovering8.m} first verifies that there is a
PQF, called $Q_{\texttt{interior}}$, which strictly satisfies the
inequalities (\ref{eq:regulator3}) for the triangulation given by
$I_1$ (see proof of Proposition~\ref{prop:periodic-triang}). By this
we know that the triangulation is a Delone triangulation. Then the
program checks that $\tilde{Q}$ satisfies all inequalities
(\ref{eq:regulator1}), (\ref{eq:regulator2}), (\ref{eq:regulator3})
given by regulators. This verifies that the Delone triangulation is in
fact a refinement of $\tilde{Q}$'s Delone subdivision.  Finally we
compute the circumradii of a representative system of simplices with
respect to the inner product $(\cdot, \cdot)$ induced by
$\tilde{Q}$. The squared circumradius $r^2$ of the simplex $L =
\conv\{\vec{0}, \vec{v}_1, \ldots, \vec{v}_8\}$ is
$$
r^2 = 
-\frac{1}{4} 
\cdot
\frac
{
\begin{vmatrix}
0 & (\vec{v}_1, \vec{v}_1) & (\vec{v}_2, \vec{v}_2) & \ldots &
(\vec{v}_8, \vec{v}_8)\\
(\vec{v}_1, \vec{v}_1) & (\vec{v}_1, \vec{v}_1) & (\vec{v}_1, \vec{v}_2) & \ldots &
(\vec{v}_1, \vec{v}_8)\\
\vdots & \vdots & \vdots & \ddots & \vdots\\
(\vec{v}_8, \vec{v}_8) & (\vec{v}_8, \vec{v}_1) & (\vec{v}_8, \vec{v}_2) & \ldots &
(\vec{v}_8, \vec{v}_8)\\
\end{vmatrix}
}
{\det\left((\vec{v}_i, \vec{v}_j)\right)_{1 \leq i,j \leq 8}},
$$
see e.g.\ \cite{sv-2004}. All these evaluations involve only rational
arithmetic and they can be carried out on a usual personal computer in
less than 15 minutes.

\end{appendix}

\bibliographystyle{amsalpha}

\vspace{2ex}

\begin{samepage}
\noindent
\textit{Achill Sch\"urmann, Department of Mathematics, University of
  Magdeburg, 39106 Magdeburg, Germany, email:
  \texttt{achill@math.uni-magdeburg.de}}
\end{samepage}

\vspace{1ex}

\begin{samepage}
\noindent
\textit{Frank Vallentin, Einstein Institute of Mathematics, The Hebrew
  University of Jerusalem, Jerusalem, 91904, Israel, email:
  \texttt{vallenti@ma.tum.de}}
\end{samepage}

\end{document}